\documentclass[11pt,a4paper,leqno]{amsart}

\usepackage{amssymb}

\usepackage{times}

\usepackage{mathrsfs}

\usepackage{verbatim}

\input{colordvi}

\theoremstyle{plain}

\newtheorem{theorem}{Theorem}[section]

\theoremstyle{remark}

\newtheorem{remark}[theorem]{Remark}

\theoremstyle{plain}

\newtheorem{lemma}[theorem]{Lemma}

\newtheorem{proposition}[theorem]{Proposition}

\numberwithin{equation}{section}

\newcommand{\al}{\alpha}

\def\Rnu{{\mathbb R}}

\def\Rn{{\Rnu^n~}}

\def\suml{\sum\limits}

\def\intl{\int\limits}

\begin{document}

\title{Relationship between stochastic flows and connection forms}

\author{Misha Neklyudov}

\date{\today}

\begin{abstract}
In this article I will prove new representation for the
Levi-Civita connection in terms of the stochastic flow
corresponding to Brownian motion on manifold.
\end{abstract}
\maketitle

The idea of reconstructing of geometry of riemannian manifold $M$
from the Brownian motion on $M$ has been productively explored for
a long time (see, for instance, expository article of Pinsky
\cite{[Pinsky-1992]} and references therein). In
\cite{[Pinsky-1992]} different asymptotics of Brownian motion
(mean exit time, distribution of exit time from small ball,...)
has been produced and it was shown that it is possible to retrieve
geometry of the manifold through the asymptotics in low dimensions
(generally less than six). Another possibility to deduce geometry
of the manifold is through small time asymptotics of logarithm of
transition function of Brownian motion on manifold (see
\cite{[Varadhan-1967]}). Indeed, Bismut-type formula (see, for
example, Corollary 3.2, p. 264 of \cite{[Elworthy-Li-1994]}) allow
us to deduce logarithm of the transition function $p:[0,T]\times
M\times M\to \mathbb{R}$ of the Brownian motion on the manifold.
We have
$$
d \log
p_t(x,y)(v_0)=\frac{1}{t}\mathbb{E}\left\{\intl_0^t<T\xi_s(v_0),X(x_s)dB_s>|\xi_t(x)=y\right\}
$$
where $\xi:[0,T]\times\Omega\to Diff(M)$ is a stochastic flow of
diffeomorphisms generated by the Brownian motion on the manifold
$$
dx_s=X(x_s)\circ dB_s+A(x_s)ds,
$$
$T\xi_{\cdot}$ is a derivative of flow $\xi_{\cdot}$ w.r.t.
initial condition. Now we can notice that $X(x_{\cdot})dB_{\cdot}$
is a martingale part of the flow $\xi$ and, consequently, can be
easily calculated. It is enough to subtract the drift of the flow
$\xi$ which, in its turn, can be calculated using Nelson
derivative of the flow.

Another way to deduce connection is arising from the theory of
Stochastic flows. It is well known fact that every nondegenerate
stochastic flow induces certain connection called Le Jan-Watanabe
connection (see example B, section 1.2 in the book
\cite{[Elworthy-Le_Jan-Li-1999]}). Indeed, connection is defined
by formulas 1.2.4 and 1.2.2 in \cite{[Elworthy-Le_Jan-Li-1999]}.
This connection is a metric connection with respect to the
riemannian metric induced by the SDE. It coincides with
Levi-Civita connection in the case of gradient Brownian systems.

My contribution is a new representation for the Levi-Civita
connection in terms of the stochastic flow, corresponding to the
Brownian motion on the manifold. Let
$\intl_{\gamma}x^idx^j,i,j=1,\ldots,n$ are areas of projections of
smooth curves $\gamma\subset M\subset\Rn$ on planes spanned by two
vectors of orthogonal basis of $\Rn$. We consider small time
asymptotic behaviour of area
$\intl_{Y_t(\gamma)}x^idx^j,i,j=1,\ldots,n$ (where $Y_t: M\to M$
is a stochastic flow of diffeomorphisms (a.s.) generated by
Brownian motion on $M$) and find connection form $\Gamma$ of
compact manifold $M$ through this asymptotic. We would like to
mention that, contrary to the Pinsky paper, our method does not
depend on dimension of manifold and it is local (i.e. knowledge of
stochastic flow in the infinitesimal neighborhood of the point $x$
immediately allow us to retrieve the connection form in $x$).

\section{Definitions and Presentation of the main result}

Let $(M,g)$ be compact riemannian manifold of dimension $k$ and
assume that it is embedded in $\Rn$, $T_xM,x\in M$ be tangent
space in the point $x\in M$, $\Gamma:\Rn\supset M\to GL(n)$,
$\Gamma=\{\Gamma_j^i(x)\}_{i,j=1}^n=\{\Gamma_{jl}^i(x)dx_l\}_{i,j=1}^n$
be Levi-Civita connection form of manifold $M$
\footnote{$\Gamma_{jl}^i$, $i,j,l=1,\ldots,n$ are Christoffel
symbols of our connection}, $P(x):\Rnu^n\to T_xM,x\in M$ be an
orthogonal projection to the tangent space. We denote
$P^{ij}(x)=(P(x)\vec{e}_i,\vec{e}_j),i,j=1,\ldots,n$, where
$\{e_i\}_{i=1}^n$ is a standard orthonormal basis in $\Rn$;
\begin{eqnarray}
S_{jl}^i(x) &=& \suml_{m=1}^nP^{im}\frac{\partial P^{jm}}{\partial
x_l}, \label{eqn:SDef}\\
r^i(x) &=& \frac{1}{2}\suml_{l=1}^nS_{il}^l,i,j,l=1,\ldots,n,x\in
M.\label{eqn:rDef}
\end{eqnarray}

Let $(\Omega,\mathcal{F},\{\mathcal{F}_t\}_{t\geq 0},\mathbb{P})$
be complete probability space with right continuous filtration,
$\{W_t\}_{t\geq 0}$ be standard Wiener process in $\Rn$. Then a
Brownian motion on the manifold $M$ is a stochastic process  $Y_t$
which satisfies following equation:
\begin{equation}
dY_t(x)=P(Y_t(x))\circ dW_t,Y_0(x)=x,x\in
M,t\in\Rnu^+,\label{eqn:BrownianMotionManifold-1}
\end{equation}
where equation is understood in Stratonovich sense, see e.g.
\cite{Stroock}, \cite{Elworthy}. We will use the same letter i.e.
$Y_t$ for the stochastic flow of diffeomorphisms corresponding to
the Brownian motion. In the notation introduced above $Y_t(A),
A\subset M$ denotes image of the set $A\subset M$ by
diffeomorphism $Y_t$.
\begin{theorem}\label{thm:MainTheorem}
Let us denote
\begin{equation}
q^i(x)=[\frac{d}{dt}\mathbb{E}Y_t^i(x)]|_{t=0},x\in
M,i=1,\ldots,n,\label{eqn:qDef}
\end{equation}
$$
\Psi^{ij}(t,\gamma)=\mathbb{E}\intl_{Y_t(\gamma)}x^idx^j,i,j=1,\ldots,
n,t\geq 0,\gamma\in C^1([0,1],M)
$$
Then $\Psi^{ij}(\cdot,\gamma),\gamma\in C^1([0,1],M)$ is
differentiable and we have following formula:
\begin{eqnarray}
\intl_{\gamma}\Gamma_{jk}^i(x)dx_k=\frac{\partial\Psi^{ij}}{\partial
t}(0,\gamma)-\frac{\partial\Psi^{ji}}{\partial
t}(0,\gamma)-2\intl_{\gamma}(q^idx^j-q^jdx^i)\nonumber\\
-(\gamma^i(1)q^j(\gamma(1))-\gamma^i(0)q^j(\gamma(0)))+(\gamma^j(1)q^i(\gamma(1))-\gamma^j(0)q^i(\gamma(0))),\nonumber\\
i,j=1,\ldots,k,\gamma\in C^1([0,1],M).\label{eqn:MainResult}
\end{eqnarray}
\end{theorem}
\begin{remark}\label{rem:DriftTerm}
It will be shown below that function $\vec{q}=(q^1,\ldots,q^n)$
can also be written as follows:
\begin{equation}
q^i(x)=\frac{1}{4}[\suml_{l=1}^n\Gamma_{il}^l(x)+\suml_{l=1}^n\frac{\partial
P^{il}}{\partial x_l}],x\in
M,i=1,\ldots,n.\label{eqn:DriftTermEstimate}
\end{equation}
\end{remark}
\begin{remark}
The formula \eqref{eqn:MainResult} allow us to find the value of
Christofell symbols $\Gamma_{jl}^i$, $i,j,l=1,\ldots,n$. 
Indeed, if $\gamma$ is a closed loop we can apply Stokes Theorem
to the left part of equality \eqref{eqn:MainResult}, divide the
result on the area of the surface and tend the size of the surface
to $0$. As the result we get the Levi-Civita connection form up to
the exact form. The remaining exact form can be calculated by
considering of curves with fixed initial point and varying end
point.
\end{remark}
The main tool for the proof of Theorem \ref{thm:MainTheorem} will
be the following proposition proved in \cite{Nekl}\footnote{We
would like to note that in \cite{Nekl} only the special case of
closed loops (i.e. $a=b$) has been considered. The general case
considered here is proved similarly.}, theorem $4$, p. 115:
\begin{proposition}\label{prop:ContourConservation}
Let $\sigma(t,\cdot)\in C_{b}^{2,\al}(\Rn ,\Rn\otimes\Rnu^{m})$,
$u(t,\cdot)\in C_{b}^{1,\al}(\Rn ,\Rn), t\in[0,T]$. Assume that
$\gamma_{a,b}$ is a curve of $C^1$ class in $\Rn$ which connects
points $a\in \Rn$ and $b\in \Rn$. Let $F\in
C^{1,2}([0,T]\times\Rn,\Rn)$,
$X=X_t(x,\omega):[0,T]\times\Rn\times\Omega\to\Rn$-be defined by:
$$
\begin{array}{l}
dX_t(x) = u(t,X_t(x))dt+\sigma(t,X_t(x))dW_t\\
X_0(x) = x.
\end{array}
$$
Then
\begin{eqnarray}\label{eqn:ContourItoFormula}
\intl_{X_t(\gamma_{a,b})}\suml_{k=1}^{n}F^k(t,x)dx_k=\intl_{\gamma_{a,b}}\suml_{k=1}^{n}F^k(0,x)dx_k+
\intl_0^t\intl_{X_s(\gamma_{a,b})}\suml_{k=1}^{n}\left(\frac{\partial
F^k}{\partial
t}\right.+\nonumber\\
\left.\suml_{j=1}^{n}u^{j}(\frac{\partial F^k}{\partial
x_j}-\frac{\partial F^j}{\partial
x_k})+\frac{1}{2}\suml_{i,j=1}^n\frac{\partial^2 F^k}{\partial
x_i\partial
x_j}\suml_{m=1}^{n}\sigma^{im}\sigma^{jm}\right)dx_kds\\
+\intl_0^t\suml_{k=1}^{n}(F^k(s,X_s(b))u^k(s,X_s(b))-F^k(s,X_s(a))u^k(s,X_s(a)))ds\nonumber\\
+\frac{1}{2}\intl_0^t\intl_{X_s(\gamma_{a,b})}\suml_{k=1}^{n}
\left(\suml_{j,l}\frac{\partial F^j}{\partial
x_l}\suml_m\sigma^{lm}\frac{\partial \sigma^{jm}}{\partial
x_k}\right)dx_kds+\intl_0^t\intl_{X_s(\gamma_{a,b})}\nonumber\\
\suml_{k,j=1}^{n}F^j(s,x)\frac{\partial \sigma^{jl}}{\partial
x_k}dx_kdw_s^l+
\intl_0^t\intl_{X_s(\gamma_{a,b})}\suml_{k=1}^n\left(\suml_{i,l=1}^n\frac{\partial
F^k}{\partial x_i}\sigma^{il}\right)dx_kdW_s^{l}.\nonumber
\end{eqnarray}
\end{proposition}
We also need
\begin{lemma}\label{lem:lem-1}
\begin{equation}
\Gamma_{jk}^i(x) = S_{jk}^i(x)-S_{ik}^j(x),i,j,k=1,\ldots,n,x\in
M\label{eqn:formula-5}
\end{equation}
\end{lemma}
\begin{proof}[Proof of lemma \ref{lem:lem-1}]
Define $S:\Rn\supset M\to GL(n)$ as follows
\begin{equation}
S(x)=\{S_j^i(x)\}_{i,j=1}^n=\{\suml_{k=1}S_{jk}^i(x)dx_k\}_{i,j=1}^n,x\in
M.\label{eqn:formula-6}
\end{equation}
We have by \eqref{eqn:SDef} that
$$
S(x)=P(x)dP^*(x).
$$
In the same time, we have (\cite{Driver},formula 3.65) that
\begin{eqnarray}
\Gamma=dQP+dPQ=-dPP+dP(Id-P)\label{eqn:formula-5'}\\
=dP-dPP-dPP=d(P^2)-2dPP=PdP-dPP=S-S^*,\nonumber
\end{eqnarray}
where  $Q=Id-P$ and we have used that $P$ is orthogonal projection
(i.e. $P^*=P$, $P^2=P$).
\end{proof}
\section{Proof of the Theorem \ref{thm:MainTheorem}}
\begin{proof}[]

We apply Proposition \ref{prop:ContourConservation} with
$F^k(t,x)=x^i\delta_{jk}$ where $i,j=1,\ldots,n$ and get
\begin{eqnarray}
\intl_{X_t(\gamma_{a,b})}x^idx_j=\intl_{\gamma_{a,b}}x^idx_j
+\intl_0^t\intl_{X_s(\gamma_{a,b})}(u^idx_j-u^jdx_i)ds\label{eqn:formula-1}\\
+\intl_0^t(X_s^i(b)u^j(s,X_s(b))-X_s^i(a)u^j(s,X_s(a)))ds+\intl_0^t\intl_{X_s(\gamma_{a,b})}\suml_{m=1}^n\sigma^{im}\frac{\partial
\sigma^{jm}}{\partial x_k}dx_kds\nonumber\\
+\intl_0^t\intl_{X_s(\gamma_{a,b})}\suml_{k,l=1}^nx^i\frac{\partial
\sigma^{jl}}{\partial
x_k}dx_kdW_s^{l}+\intl_0^t\suml_{l=1}^n\intl_{X_s(\gamma_{a,b})}\sigma^{il}dx_jdW_s^{l}.\nonumber
\end{eqnarray}
Taking mathematical expectation of both parts of formula
\eqref{eqn:formula-1} we get
\begin{eqnarray}
\mathbb{E}\intl_{X_t(\gamma_{a,b})}x^idx_j=\intl_{\gamma_{a,b}}x^idx_j
+\mathbb{E}\intl_0^t\intl_{X_s(\gamma_{a,b})}(u^idx_j-u^jdx_i)ds\nonumber\\
+\mathbb{E}\intl_0^t(X_s^i(b)u^j(s,X_s(b))-X_s^i(a)u^j(s,X_s(a)))ds\nonumber\\
+\mathbb{E}\intl_0^t\intl_{X_s(\gamma_{a,b})}\suml_{m=1}^n\sigma^{im}\frac{\partial
\sigma^{jm}}{\partial x_k}dx_kds.\label{eqn:formula-2}
\end{eqnarray}
Let us rewrite equation \eqref{eqn:BrownianMotionManifold-1} for
Brownian motion $\{Y_t\}_{t\geq 0}$ on $M$ in the Ito form. We
have
\begin{eqnarray}
dY_t^i(x)=\suml_{j=1}^nP^{ij}(Y_t(x))\circ
dW_t^j\nonumber\\
=\suml_{j=1}^nP^{ij}(Y_t(x))dW_t^j+\frac{1}{2}\suml_{l,j=1}^nP^{lj}\frac{\partial
P^{ij}}{\partial
x_l}(Y_t(x))dt\nonumber\\
=\suml_{j=1}^nP^{ij}(Y_t(x))dW_t^j+\frac{1}{2}\suml_{l=1}^nS_{il}^l(Y_t(x))dt\nonumber\\
=r^i(Y_t(x))dt+\suml_{j=1}^nP^{ij}(Y_t(x))dW_t^j,i=1,\ldots,n,t\geq
0.\label{eqn:BrownianMotionManifold-2}
\end{eqnarray}
Now we can put $X_t=Y_t$ in the formula \eqref{eqn:formula-2} and
we get
\begin{eqnarray}
\Psi^{ij}(t,\gamma_{a,b})=\intl_{\gamma_{a,b}}x^idx_j
+\mathbb{E}\intl_0^t\intl_{Y_s(\gamma_{a,b})}(r^idx_j-r^jdx_i)ds\label{eqn:formula-3}\\
+\mathbb{E}\intl_0^t(Y_s^i(b)r^j(Y_s(b))-Y_s^i(a)r^j(Y_s(a)))ds
+\mathbb{E}\intl_0^t\intl_{Y_s(\gamma_{a,b})}S_{jk}^i(x)dx_kds.\nonumber
\end{eqnarray}
Therefore, we have
\begin{eqnarray}
\frac{\partial\Psi^{ij}}{\partial
t}(0,\gamma_{a,b})=\intl_{\gamma_{a,b}}(r^idx_j-r^jdx_i)+(b^ir^j(b)-a^ir^j(a))\nonumber\\
+\intl_{\gamma_{a,b}}S_{jk}^i(x)dx_k,i,j=1,\ldots,n.\label{eqn:formula-4}
\end{eqnarray}

It remains to show that
$$
r^i=q^i,i=1,\ldots,n
$$
and formula \eqref{eqn:MainResult} will immediately follow from
\eqref{eqn:formula-4} and \eqref{eqn:formula-5}. We can notice
that $r$ is a drift of Brownian motion by formula
\eqref{eqn:BrownianMotionManifold-2}. Now applying mathematical
expectation to formula \eqref{eqn:BrownianMotionManifold-2} we
immediately get the result.
\end{proof}
\begin{proof}[Proof of Remark \ref{rem:DriftTerm}]
We have
\begin{equation*}
S+S^*=PdP+dPP=d(P^2)=dP
\end{equation*}
i.e.
\begin{equation}
S_{jm}^i(x)+S_{im}^j(x)=\frac{\partial P^{ij}}{\partial
x_m},ij,m=1,\ldots,n.\label{eqn:formula-7}
\end{equation}
Therefore,from \eqref{eqn:formula-5} and \eqref{eqn:formula-7} it
follows that
\begin{eqnarray}
r^i=\frac{1}{2}\suml_{l=1}^nS_{il}^l=\frac{1}{4}\suml_{l=1}^n[(S_{il}^l-S_{ll}^i)+(S_{il}^l-S_{ll}^i)]\nonumber\\
=\frac{1}{4}\suml_{l=1}^n[\Gamma_{il}^l+\frac{\partial
P^{il}}{\partial x_l}]=q_i,i=1,\ldots,n.
\end{eqnarray}
\end{proof}

\section{Acknowledgement}

The Author would like to thank D. Elworthy and Z. Brzezniak for
making valuable suggestions.

\end{document}